\documentclass[]{spie}  

\usepackage[utf8]{inputenc}
\usepackage[english]{babel}

\usepackage{epsfig}
\usepackage{ellipsis}
\usepackage{amsmath}
\usepackage{amssymb}
\usepackage{url}
\usepackage{array}
\usepackage{color}
\usepackage{multirow}
\usepackage{epsfig}
\usepackage{colortbl}
\usepackage [table]{xcolor}
\usepackage{graphicx}
\usepackage{subcaption}
\usepackage{float}
\usepackage{caption}
\usepackage{algorithm}
\usepackage[algo2e]{algorithm2e}
\usepackage{mathtools}
\usepackage{rotating} 
\usepackage{textcomp} 
\usepackage{nomencl} 
\makenomenclature
\SetKwData{Left}{left}
\SetKwData{This}{this}
\SetKwData{Up}{up}
\SetKwFunction{Union}{Union}
\SetKwFunction{FindCompress}{FindCompress}
\SetKwInOut{Input}{Input}
\SetKwInOut{Output}{Output}

\DeclareMathOperator*{\argmin}{\arg\!\min}

\DeclareMathOperator{\sgn}{sgn}

\usepackage[]{graphicx}
\usepackage{subcaption}
\usepackage{epstopdf}
\usepackage{cite}
\usepackage{textcomp}
\usepackage{cite}
\usepackage{amsmath}
\usepackage{amsfonts}
\usepackage{mathtools}
\usepackage{algorithm}
\usepackage{algpseudocode}
\usepackage[utf8]{inputenc}
\usepackage[english]{babel}
\graphicspath{ {figures/}}

\title{Stochastic First-Order Minimization Techniques Using Jensen Surrogates for X-Ray Transmission Tomography} 

\author{Soysal Degirmenci\supit{a}, Joseph A. O'Sullivan\supit{a}, David G. Politte\supit{b} 
\skiplinehalf
\supit{a}Department of Electrical and Systems Engineering, Washington University,\\ 1 Brookings Drive, Saint Louis, MO, USA, 63130 \\
\supit{b}Mallinckrodt Institute of Radiology, Washington University School of Medicine,\\ 510 South Kingshighway Blvd, Saint Louis, MO, USA, 63110
}

\authorinfo{\\Soysal Degirmenci: E-mail: s.degirmenci@wustl.edu\\Joseph A. O'Sullivan: E-mail: jao@wustl.edu\\David G. Politte: E-mail: politted@wustl.edu}

\begin{document} 
\maketitle 

\begin{abstract}
Image reconstruction in X-ray transmission tomography has been an important research field for decades. In light of data volume increasing faster than processor speeds, one needs accelerated iterative algorithms to solve the optimization problem in the X-ray CT application. Incremental methods, in which a subset of data is being used at each iteration to accelerate the computations, have been getting more popular lately in the machine learning and mathematical optimization fields. The most popular member of this family of algorithms in the X-ray CT field is ordered-subsets. Even though it performs well in earlier iterations, the lack of convergence in later iterations is a known phenomenon. In this paper, we propose two incremental methods that use Jensen surrogates for the X-ray CT application, one stochastic and one ordered-subsets type. Using measured data, we show that the stochastic variant we propose outperforms other algorithms, including the gradient descent counterparts.
\end{abstract}


\section{Introduction}

With the scale of computational inverse problems getting larger in many applications, including X-Ray Computed Tomography, acceleration techniques that provide fast rates of convergence have become crucial. One of the acceleration techniques widely used in many applications is incremental methods, which is a range decomposition technique~\cite{bertsekas2011incremental}. In the general case, one forms a surrogate function around the current estimate using a subset of data and computes the next estimate by minimizing it. The idea behind the acceleration is that the subset of data is a good approximation of the full data space and computational time is saved by using only a portion of it. The most popular variants of incremental methods, ordered subsets, and stochastic surrogate descent, perform well in earlier iterations but are known to lack convergence for the later iterations. Lately, the variants of incremental methods in optimization field have proliferated. The algorithms proposed were shown to perform well with desirable convergence properties. The idea behind these methods is to keep track of the surrogate parameters for each subset and to use a combination of them while performing updates, rather than using only one of them.

Jensen surrogates is a possible choice of surrogate functions, and it has nice properties that result in algorithms with attractive properties and good performance~\cite{o2007alternating, degirmenci2015acceleration}. For X-ray CT, the family of algorithms that uses Jensen surrogates has parameters that are easy to compute. These can also be extended to acceleration methods, including range-based techniques.

In this paper, we:
\begin{itemize}
	\item explore stochastic incremental methods and propose two novel accelerated algorithms using Jensen surrogates for X-ray CT, Stochastic Average Jensen Surrogates Optimization (SA-JS) and Ordered Subsets Average Jensen Surrogates Optimization (OSA-JS), and
	\item show that the SA-JS algorithm outperforms the other competing algorithms for a large number of subsets by using two sets of real data collected from a baggage scanner.
\end{itemize}

\section{Convex Optimization Using Jensen Surrogates for X-Ray CT}
 The objective function we minimize for monoenergetic X-Ray transmission tomography is as follows.
\begin{eqnarray} \label{eq8.1.1.1}
\min_{\boldsymbol{x} \in \mathbb{X}} \Phi(\boldsymbol{x}) = \min_{\boldsymbol{x} \in \mathbb{X}} f(\boldsymbol{x}) + \lambda \beta(\boldsymbol{x}) &=& \min_{\boldsymbol{x} \in \mathbb{X}} \sum_{i=1}^{M}f_{i}(\boldsymbol{x}) +  \lambda \sum_{k=1}^{K} \beta_{k}(\boldsymbol{x}) \\ 
&=& \min_{\boldsymbol{x} \in \mathbb{R}_{+}^{N}} \sum_{i=1}^{M}\tilde{f}_{i}((\boldsymbol{H}\boldsymbol{x})_{i}) +  \lambda \sum_{k=1}^{K} \tilde{\beta}_{k}((\boldsymbol{C}\boldsymbol{x})_{k}),
\end{eqnarray}
where $f(\boldsymbol{x})$ is the data-fitting term and $\beta(\boldsymbol{x})$ is the regularization term. Here, we specifically look at the Poisson log-likelihood data-fitting term, where
\begin{eqnarray} \label{eq8.1.1.2}
\tilde{f}_{i}(l) = d_{i} l + I_{0,i} \exp(-l),
\end{eqnarray}
$\boldsymbol{d} \in \mathbb{R}_{+}^{M}$ is the attenuated data vector, $\boldsymbol{I_{0}} \in \mathbb{R}_{+}^{M}$ is the incident photon count vector, and $\boldsymbol{H} \in \mathbb{R}_{+}^{M\times N} $ is the system matrix that defines the relationship between ray-paths and voxels. In a simple ray-tracing model, $h_{ij}$ represents the length of intersection between the voxel indexed by $j$ and the ray-path indexed by $i$. The non-negativity constraint on the image is due to the physical nature of linear attenuation coefficients. The regularization term we use is
\begin{eqnarray} \label{eq8.1.1.3}
\tilde{\beta}_{k}(t) = \omega_{k} \delta^{2} \Bigg( \Big| \frac{t}{\delta} \Big| - \log (1 + \Big| \frac{t}{\delta} \Big| ) \Bigg),
\end{eqnarray}
where $\omega_{k} > 0$, $\delta > 0$, and $\boldsymbol{C} \in \mathbb{R}^{K \times N}$ is a matrix that has $1$s along the diagonal and only one $-1$ off-diagonal for each row. We further assume that for two arbitrary voxels, if there exists a row in $\boldsymbol{C}$ where the first voxel has a value equal to $1$ and the second $-1$, there should exist another row where the reverse is true as well. The $\tilde{\beta}_{k}$ defined is an edge-preserving function that is of Huber type, convex, even, and differentiable, where $\boldsymbol{C}$ defines a neighborhood around a center voxel we would like to use to regularize the image. Similarly, $\omega_{k}$ determines the weights of the corresponding neighborhood. This function behaves quadratically for small $|t/\delta|$ and linearly for large $|t/\delta|$ values. Since both $\tilde{f}_{i}$ and $\tilde{\beta}_{k}$  are convex and differentiable, we can use the same Jensen surrogate formulation for each of them to create an iterative minimization algorithm.

For the data-fitting term, $f(\boldsymbol{x})$, we denote its Jensen surrogate around $\boldsymbol{x}$ parameterized by $\boldsymbol{r}$ as $g_{\boldsymbol{r}}(\boldsymbol{x}; \boldsymbol{\hat{x}})$. We follow the same procedure as in,\cite{degirmenci2015acceleration} which gives the resulting surrogate\footnote{We ignore the constant term after step 1.}

\begin{eqnarray} \label{eq8.1.1.44}
g_{\boldsymbol{r}}(\boldsymbol{x}; \boldsymbol{\hat{x}}) = \sum_{i}g_{i,\boldsymbol{r}_{i}}(\boldsymbol{x}; {\boldsymbol{\hat{x}}}) &=& \sum_{i} \sum_{j}r_{ij} \tilde{f}_{i} \Big( \frac{h_{ij}}{r_{ij}} (x_{j}-\hat{x}_{j}) + (\boldsymbol{H}\boldsymbol{\hat{x}})_{i}  \Big) + \text{const.} \nonumber \\
&=& \sum_{i} \sum_{j} h_{ij} d_{i} (x_{j}-\hat{x}_{j}) + r_{ij} I_{0,i} \exp(-\frac{h_{ij}}{r_{ij}} (x_{j}-\hat{x}_{j}) - (\boldsymbol{H}\boldsymbol{\hat{x}})_{i}).
\end{eqnarray}
We choose $r_{ij} = h_{ij}/Z$, where $Z = \max_{i} \sum_{j} h_{ij}$, which results in a surrogate function as follows.
\begin{eqnarray} \label{eq8.1.1.4}
g_{\boldsymbol{r}}(\boldsymbol{x}; \boldsymbol{\hat{x}}) = \sum_{i}g_{i,\boldsymbol{r}_{i}}(\boldsymbol{x}; {\boldsymbol{\hat{x}}}) &=& \sum_{i} \sum_{j}\frac{h_{ij}}{Z} \tilde{f}_{i} \Big( Z(x_{j}-\hat{x}_{j}) + (\boldsymbol{H}\boldsymbol{\hat{x}})_{i}  \Big) + \text{const.} \nonumber \\
&=& \sum_{i} \sum_{j} h_{ij} d_{i} (x_{j}-\hat{x}_{j}) + \frac{h_{ij}}{Z} I_{0,i} \exp(-Z(x_{j}-\hat{x}_{j}) - (\boldsymbol{H}\boldsymbol{\hat{x}})_{i}) \nonumber \\
&=& \sum_{i} \sum_{j} h_{ij} d_{i} (x_{j}-\hat{x}_{j}) + \frac{h_{ij}}{Z} \hat{q}_{i} \exp(-Z(x_{j}-\hat{x}_{j})) \nonumber \\
&=& \sum_{i} b_{j} (x_{j}-\hat{x}_{j}) + \frac{\hat{b}_{j}}{Z} \exp(-Z(x_{j}-\hat{x}_{j})).
\end{eqnarray}

Here, it is important to note that when $r_{ij} = \frac{h_{ij}\hat{x}_{j}}{\sum_{j'}h_{ij'}\hat{x}_{j'}}$, the surrogate does not have a closed-form minimization and thus must be solved using some convex minimization method. Lange, Fessler et al.~\cite{lange1995globally} use this auxiliary variable choice and solves it via Newton's method.

For the regularization term $\beta(\boldsymbol{x})$ we denote its Jensen surrogate around $\boldsymbol{x}$ parameterized by $\boldsymbol{s}$ as $B_{\boldsymbol{s}}(\boldsymbol{x}; \boldsymbol{\hat{x}})$. We choose $s_{ij} = |c_{ij}|/2$ (since $c_{ij} \in \{-1, 0, 1 \}$ and there is exactly one element equal to $1$ and one element equal to $-1$ in each row of $\boldsymbol{C}$, the denominator is equal to 2). Then the surrogate becomes\footnote{We ignore the constant term after step 1.}

\begin{eqnarray} \label{eq8.1.1.5}
B_{\boldsymbol{s}}(\boldsymbol{x}; \boldsymbol{\hat{x}}) = \sum_{k}B_{k,\boldsymbol{s}}(\boldsymbol{x}; {\boldsymbol{\hat{x}}}) &=& \sum_{k} \sum_{j}s_{kj} \tilde{\beta}_{k} \Big( \frac{c_{kj}}{s_{kj}}(x_{j}-\hat{x}_{j}) + (\boldsymbol{C}\boldsymbol{\hat{x}})_{k}  \Big) + \text{const.} \nonumber \\
&=& \sum_{k} \sum_{j} \frac{|c_{kj}|}{2} \tilde{\beta}_{k} \Big( 2 \sgn{(c_{kj})} (x_{j}-\hat{x}_{j}) + (\boldsymbol{C}\boldsymbol{\hat{x}})_{k}  \Big) \nonumber \\
&=& \sum_{j' \in N_{j}} \sum_{j}  \frac{|c_{j'j}|}{2} \tilde{\beta}_{j'} \Big( 2 \sgn{(c_{j'j})} (x_{j}-\hat{x}_{j}) + \sgn{(c_{j'j})} (\hat{x}_{j} -\hat{x}_{j'})  \Big) \nonumber \\
&=& \sum_{j} \sum_{j' \in N_{j}} \frac{|c_{j'j}|}{2} \tilde{\beta}_{jj'} \Big( 2 \sgn{(c_{j'j})} (x_{j}-\hat{x}_{j}) + \sgn{(c_{j'j})} (\hat{x}_{j} -\hat{x}_{j'})  \Big) \nonumber \\
&=& \sum_{j} \sum_{j' \in N_{j}} \frac{1}{2} \tilde{\beta}_{jj'} \Big( 2 (x_{j}-\hat{x}_{j}) +  (\hat{x}_{j} -\hat{x}_{j'})  \Big) \nonumber \\
&=& \sum_{j} \sum_{j' \in N_{j}} \frac{1}{2} \tilde{\beta}_{jj'} \Big( 2 x_{j} - \hat{x}_{j} -\hat{x}_{j'}  \Big) \nonumber \\
&=& \sum_{j} \sum_{j' \in N_{j}} \frac{\omega_{jj'}}{2} \delta^{2} \Bigg( \Big| \frac{2 x_{j} - \hat{x}_{j} -\hat{x}_{j'}}{\delta} \Big| - \nonumber \\  
& & \log (1 + \Big| \frac{2 x_{j} - \hat{x}_{j} -\hat{x}_{j'}}{\delta} \Big| ) \Bigg),
\end{eqnarray}
where in the third step we changed the notation, so that $N_{j}$ is a set of indices that defines the neighborhood of index $j$ (those sets of indices are the only ones that are relevant to $x_{j}$.), and in the fifth step we used the fact that $\tilde{\beta}$ is an even function.

When two surrogate functions are combined, we have our algorithm, the Jensen Surrogates Optimization for X-Ray CT, presented in Algorithm~\ref{algo_full}. 
\begin{algorithm} 
	\caption{Jensen Surrogates Optimization for X-Ray CT (Full-JS)}
	\label{algo_full}
	\Input{$\boldsymbol{x}^{(0)} \in \mathbb{R}_{+}^{N}, \boldsymbol{d}, \boldsymbol{I_{0}} \in \mathbb{R}_{+}^{M}, \boldsymbol{H} \in \mathbb{R}^{M \times N}, \lambda \ge 0, \delta > 0$.}
	\textbf{Pre-compute} $b_{j} = \sum_{i} d_{i}h_{ij}$, $\forall j$. \\
	\textbf{Pre-compute} $Z = \max_{i} \sum_{j} h_{ij}$. \\
	\For {$n = 0, 1, 2, ...$}{
		$q_{i}^{(n)} = I_{0,i} \exp(-\sum_{j}h_{ij}x^{(n)}_{j})$, $\forall i$. \\
		${b}^{(n)}_{j} = \sum_{i} q_{i}^{(n)} h_{ij}$, $\forall j$. \\
		$x_{j}^{(n+1)} = \argmin_{x \ge 0} b_{j} (x - x^{(n)}_{j}) + {b}^{(n)}_{j}/Z \exp(-Z(x -x^{(n)}_{j})) + \lambda \sum_{j' \in N_{j}} \frac{1}{2} \tilde{\beta}_{jj'} ( 2 x - \hat{x}_{j} -\hat{x}_{j'})$, $\forall j$.
	}
\end{algorithm}

As in most iterative algorithms, this algorithm consists of three steps in each iteration. First, the current estimate in the image domain is forward projected to the data domain with Beer's Law of attenuation applied. Then, this estimate is back projected onto the image domain where this information as well as the back projection of attenuated data and the previous image estimate are used to compute the next estimate. In the case of no regularization ($\lambda = 0$), the surrogate can be minimized in one step; there is a closed-form update. For the regularized case, however, there is not a closed-form update that minimizes the combined surrogate functions. Instead, we have $N$ independent one-dimensional convex problems we can minimize in parallel. Any convex optimization method can be used to minimize these functions. In order to achieve fast convergence we first attempted to use Newton's method. It is important to note that since each minimization is a one-dimensional problem, inversion of the Hessian is not an issue and we take advantage of that fact. However, due to characteristics of the functions being minimized, Newton's method diverged for some cases. Then, we attempted to use a Trust Region Newton's method,\cite{nocedal2006numerical} which is a modification of Newton's method. In the trust region method, in each iteration, there is a trust region defined such that it bounds the next iterate computed, which makes the algorithm more stable. A metric that measures how well the quadratic approximation approximates the original function is computed at each iteration. Depending on the value of this metric, we either ``trust" the quadratic approximation more and expand the trust region, or trust less and shrink it.  Unfortunately, the trust region method requires two computations of the function, one first derivative, one second derivative, and many comparisons per iteration. In order to get a faster method that requires fewer computations per iteration and produces comparable performance, we developed a modified trust region method. This method takes advantage of the structure of the problem to construct a fixed trust region and only requires one first derivative, one second derivative, and two comparisons per iteration. More information about these methods can be found in Degirmenci.\cite{degirmenci2016thesis}

As discussed above, when $\lambda = 0$, the image update at iteration $n$ becomes
\begin{eqnarray} \label{eq8.1.1.6}
x_{j}^{(n+1)} = \Bigg[ x_{j}^{(n)} - \frac{1}{Z} \log \Big( \frac{b_{j}}{{b}^{(n)}_{j}} \Big) \Bigg]_{+} \text{, } \forall j,
\end{eqnarray}
where $[\cdot]_{+} = \max(0, \cdot)$ is a non-negativity operator. Unfortunately, these types of ``full" methods have a sublinear rate of convergence~\cite{degirmenci2016thesis}, which is very slow. Applications of X-ray imaging, such as baggage scanning (considered below) and medical imaging, both require fast reconstructions with ``good" resultant image volumes. In the next section, we will look at several acceleration methods based on range decompositions, which are variants of Algorithm \ref{algo_full}.

\section{Acceleration Methods}

When the amount of data is very large, a reasonable approach to reach the optimal solution faster is to use a portion or subset of the data at each iteration to perform updates. These methods are called ``incremental methods" in general. In this section, we first present the popular incremental approach in X-ray CT called ordered subsets using Jensen surrogates and propose our new algorithms Stochastic Average and Ordered Subsets Average Jensen Surrogates Optimization for X-Ray CT. For convenience, we assume that the data is split into $B^{r}$ subsets where the set of source-detector indices of subset $k$ is represented as $\mathcal{B}^{r}_{k}$.

\subsection{Ordered Subsets Jensen Surrogates Optimization for X-Ray CT}
Ordered subsets is a popular technique used in X-ray CT in which the estimates are updated at each iteration by using a subset of data in an ordered way. In other words, the choice of indices to be used is deterministic and follows a cyclic order. Since only a subset of data is used, in order to keep the balance between the data-fitting and the regularization term, the regularization parameter is scaled down by $B^{r}$, the number of subsets. The corresponding algorithm of ordered subsets when used with Jensen surrogates is presented in Algorithm~\ref{algo_os}.
\begin{algorithm} 
	\caption{Ordered Subsets Jensen Surrogates Optimization for X-Ray CT (OS-JS)}
	
	\Input{$\boldsymbol{x}^{(0)} \in \mathbb{R}_{+}^{N}, \boldsymbol{d}, \boldsymbol{I_{0}} \in \mathbb{R}_{+}^{M}, \boldsymbol{H} \in \mathbb{R}^{M \times N}, \lambda \ge 0, \delta > 0$, $\mathcal{B}^{r}_{k}$ for $k=0, 1, ..., (B^{r}-1)$}
	\textbf{Pre-compute} $b^{k}_{j} = \sum_{i \in \mathcal{B}^{r}_{k} } d_{i}h_{ij}$, $\forall j, k$. \\
	\textbf{Pre-compute} $Z = \max_{i} \sum_{j} h_{ij}$. \\
	\For {$n = 0, 1, 2, ...$}{
		$k = \mod{(n, B^{r})}$ \\
		$q_{i}^{(n)} = I_{0,i} \exp(-\sum_{j}h_{ij}x^{(n)}_{j})$, $\forall i \in \mathcal{B}^{r}_{k} $. \\
		${b}^{(n, k)}_{j} = \sum_{i \in \mathcal{B}^{r}_{k} } q_{i}^{(n)} h_{ij}$, $\forall j$. \\
		$x_{j}^{(n+1)} = \argmin_{x \ge 0} b^{k}_{j} (x - x^{(n)}_{j}) + {b}^{(n, k)}_{j}/Z \exp(-Z(x -x^{(n)}_{j})) + \frac{\lambda}{B^{r}} B_{\boldsymbol{s}}(\boldsymbol{x}; \boldsymbol{x}^{(n)})$, $\forall j$.
	}
	\label{algo_os}
\end{algorithm}
This algorithm is known to work well in earlier iterations when the estimate is ``far away" from the optimum. However, for later iterations, it lacks convergence and stops decreasing after some number of iterations. For the gradient descent case, Bertsekas et al.~\cite{bertsekas2011incremental} showed that regardless of the number of iterations, this type of method never gets closer to the optimal function value than a certain positive constant.

\subsection{Stochastic Average Jensen Surrogates Optimization for X-Ray CT}
As seen in the previous section, when the ordered subsets technique is used, we use only the most recent version of the back projection information in order to minimize a subset of the problem. In this section, we propose to store the back projection of the estimates for each subset, whether they are the most recent or not, and to use a combination of this back projection information to update the estimates at each iteration. The subset is chosen stochastically, so that each subset has a probability of being chosen equal to $1/B^{r}$. This algorithm can be seen as an equivalent of the Stochastic Average Gradient algorithm~\cite{roux2012stochastic} that uses Jensen surrogates instead of Lipschitz quadratic surrogates for the X-Ray CT application. Also, it is important to note that we only use the previous estimate information for the data-fitting term; a similar approach can be used for the regularization term as well, but this would require the storage of $B^{r}$ images. Not only does this make the one-dimensional optimization to update the estimate slower, but also increases the storage requirement. Experimentally, we have not found any advantages and will not include it here.

Algorithm~\ref{algo_sa} presents the Stochastic Average Jensen Surrogates Optimization for X-Ray CT (SA-JS). For a very large number of subsets, the requirement of storing $B^{r}$ back projected estimates in the image domain can be problematic if the memory capacity of the computing architecture being used is low. As seen in the image update, at each iteration we use the most recent sum of the back projected estimates. In a naive implementation, this computation would add $\mathcal{O}(B^{r}N)$ to the computational complexity. A better approach is to store the most recent sum as another variable, and when a subset is chosen, subtract the old back projected estimate of the corresponding estimate, and add the newly computed one to the sum. This version would require $\mathcal{O}(2N)$ in additional computations and an additional $\mathcal{O}(N)$ of storage, which is more reasonable than the naive implementation.

\begin{algorithm} 
	\caption{Stochastic Average Jensen Surrogates Optimization for X-Ray CT (SA-JS)}
	
	\Input{$\boldsymbol{x}^{(0)} \in \mathbb{R}_{+}^{N}, \boldsymbol{d}, \boldsymbol{I_{0}} \in \mathbb{R}_{+}^{M}, \boldsymbol{H} \in \mathbb{R}^{M \times N}, \lambda \ge 0, \delta > 0$, $\mathcal{B}^{r}_{k}$ for $k=0, 1, ..., (B^{r}-1)$, $\boldsymbol{x}^{(0, k)} \in \mathbb{R}_{+}^{N}$ for $k=0, 1, ..., (B^{r}-1)$}
	\textbf{Pre-compute} $b^{k}_{j} = \sum_{i \in \mathcal{B}^{r}_{k} } d_{i}h_{ij}$, $\forall j, k$. \\
	\textbf{Pre-compute} $Z = \max_{i} \sum_{j} h_{ij}$. \\
	\For {$n = 0, 1, 2, ...$}{
		\text{Choose $k$ from $\{0, 1, ..., (B^{r}-1)\}$ randomly.}\\
		$\boldsymbol{x}^{(n, k)} = \boldsymbol{x}^{(n)}$ \\
		$q_{i}^{(n)} = I_{0,i} \exp(-\sum_{j}h_{ij}x^{(n)}_{j})$, $\forall i \in \mathcal{B}^{r}_{k} $. \\
		${b}^{(n, k)}_{j} = \sum_{i \in \mathcal{B}^{r}_{k} } q_{i}^{(n)} h_{ij}$, $\forall j$. \\
		$x_{j}^{(n+1)} = \argmin_{x \ge 0} \sum_{k} b^{k}_{j} (x - x^{(n)}_{j}) + \sum_{k} {b}^{(n, k)}_{j}/Z \exp(-Z(x -x^{(n)}_{j})) + \lambda B_{\boldsymbol{s}}(\boldsymbol{x}; \boldsymbol{x}^{(n)})$, $\forall j$.
	}
	\label{algo_sa}
\end{algorithm}
\subsection{Ordered Subsets Average Jensen Surrogates Optimization for X-Ray CT}
In addition to the stochastic choice of subsets presented in the previous section, one can also choose them in a cyclic manner as in the ordered subsets method. The resultant algorithm, called Ordered Subsets Average Jensen Surrogates Optimization for X-Ray CT (OSA-JS), is shown in Algorithm~\ref{algo_osa}.

\begin{algorithm} 
	\caption{Ordered Subsets Average Jensen Surrogates Optimization for X-Ray CT (OSA-JS)}
	
	\Input{$\boldsymbol{x}^{(0)} \in \mathbb{R}_{+}^{N}, \boldsymbol{d}, \boldsymbol{I_{0}} \in \mathbb{R}_{+}^{M}, \boldsymbol{H} \in \mathbb{R}^{M \times N}, \lambda \ge 0, \delta > 0$, $\mathcal{B}^{r}_{k}$ for $k=0, 1, ..., (B^{r}-1)$, $\boldsymbol{x}^{(0, k)} \in \mathbb{R}_{+}^{N}$ for $k=0, 1, ..., (B^{r}-1)$}
	\textbf{Pre-compute} $b^{k}_{j} = \sum_{i \in \mathcal{B}^{r}_{k} } d_{i}h_{ij}$, $\forall j, k$. \\
	\textbf{Pre-compute} $Z = \max_{i} \sum_{j} h_{ij}$. \\
	\For {$n = 0, 1, 2, ...$}{
		$k = \mod{(n, B^{r})}$ \\
		$\boldsymbol{x}^{(n, k)} = \boldsymbol{x}^{(n)}$ \\
		$q_{i}^{(n)} = I_{0,i} \exp(-\sum_{j}h_{ij}x^{(n)}_{j})$, $\forall i \in \mathcal{B}^{r}_{k} $. \\
		${b}^{(n, k)}_{j} = \sum_{i \in \mathcal{B}^{r}_{k} } q_{i}^{(n)} h_{ij}$, $\forall j$. \\
		$x_{j}^{(n+1)} = \argmin_{x \ge 0} \sum_{k} b^{k}_{j} (x - x^{(n)}_{j}) + \sum_{k} {b}^{(n, k)}_{j}/Z \exp(-Z(x -x^{(n)}_{j})) + \lambda B_{\boldsymbol{s}}(\boldsymbol{x}; \boldsymbol{x}^{(n)})$, $\forall j$.
	}
	\label{algo_osa}
\end{algorithm}

\section{Results}
In this section, we compare the proposed algorithms and the gradient descent variants using real data. Before presenting the results, we briefly explain the competing algorithms implemented. 

The gradient descent method is a popular technique for iterative convex minimization. In an iteration, one forms a surrogate function around the current estimate that is equal to
\begin{eqnarray}
	\Phi (\boldsymbol{x}^{(n)}) + \nabla (\Phi (\boldsymbol{x}^{(n)}))^{T}(\boldsymbol{x} - \boldsymbol{x}^{(n)}) + \frac{\alpha^{(n)}}{2} \| \boldsymbol{x} - \boldsymbol{x}^{(n)} \|_{2}^{2}.
\end{eqnarray}
When the function of interest is twice differentiable and Lipschitz continuous, this method has a guarantee of convergence when the quadratic scaling factor $\alpha^{(n)}$ is chosen to be equal to $L_{\Phi}$, where $L_{\Phi}$ is the Lipschitz gradient constant of $\Phi$.\footnote{Using an arbitrary initial $\alpha^{(n)}$ and decreasing it when the objective function starts increasing is a popular technique when this constant is not computable. This method, backtracking, can also be used with Jensen surrogates. For simplicity, this investigation is left as future work.} For our case, the gradient is equal to
\begin{eqnarray}
\nabla \Phi(\boldsymbol{\hat{x}}) &=& \nabla f(\boldsymbol{\hat{x}}) + \lambda \nabla \beta(\boldsymbol{\hat{x}}) \\
&=& \boldsymbol{H}^{T} (\boldsymbol{d} - \boldsymbol{\hat{q}}) + \lambda \boldsymbol{C}^{T} (\nabla \tilde{\beta}(\boldsymbol{C} \boldsymbol{\hat{x}})),
\end{eqnarray}
where $\hat{q}_{i} = I_{0,i}\exp(-\sum_{j}h_{ij}\hat{x}_{j})$. The Lipschitz gradient constant $L_{\Phi}$ is computed using the following inequality on the Hessian,
\begin{eqnarray}
\nabla^{2} \Phi(\boldsymbol{x}) &=& \nabla^{2} f(\boldsymbol{x}) + \lambda \nabla^{2} \beta(\boldsymbol{x}) \\
&\le& \max_{l, i}\frac{\partial^{2} \tilde{f}_{i}(l)}{\partial l^{2}} \boldsymbol{H}^{T}\boldsymbol{H} + \lambda \max_{l, k}\frac{\partial^{2} \tilde{\beta}_{k}(l)}{\partial l^{2}} \boldsymbol{C}^{T}\boldsymbol{C} \\
&\le& \max_{i}I_{0, i} \boldsymbol{H}^{T}\boldsymbol{H} + \lambda \boldsymbol{C}^{T}\boldsymbol{C},
\end{eqnarray}
which in return can be used to find $L_{\Phi}$ as
\begin{eqnarray} \label{eq_lip_xray}
L_{\Phi} = \Lambda_{max}(\max_{i}I_{0, i} \boldsymbol{H}^{T}\boldsymbol{H} + \lambda \boldsymbol{C}^{T}\boldsymbol{C}) \ge \max_{\boldsymbol{x}} \Lambda_{max}(\nabla^{2} \Phi(\boldsymbol{x})),
\end{eqnarray}
where $\Lambda_{max}(\cdot)$ is the operator that returns the maximum eigenvalue of the given square matrix. When the system matrix is scaled down in the image domain, we experimentally found that $\Lambda_{max}(\boldsymbol{H}^{T}\boldsymbol{H})$ follows a linear relationship with the fully sampled version. Thus, we computed it by finding $\Lambda_{max}(\boldsymbol{H}^{T}\boldsymbol{H})$ for a scaled-down system matrix and then multiplying by the appropriate scale factor. For the regularization matrix part, we computed $\Lambda_{max}(\boldsymbol{C}^{T}\boldsymbol{C})$ for a subset of the image volume, $5$ z-slices, to make it computationally tractable.

Algorithm~\ref{algo_gd_full} presents the Gradient Descent Optimization for X-Ray CT (Full-GD). For convenience, we will compare the ordered subset (OS-GD) and the stochastic average (SA-GD) variants of the gradient descent method, which are straightforward extensions of the base case algorithm.
\begin{algorithm} 
	\caption{Gradient Descent Optimization for X-Ray CT (Full-GD)}
	\label{algo_gd_full}
	\Input{$\boldsymbol{x}^{(0)} \in \mathbb{R}_{+}^{N}, \boldsymbol{d}, \boldsymbol{I_{0}} \in \mathbb{R}_{+}^{M}, \boldsymbol{H} \in \mathbb{R}^{M \times N}, \lambda \ge 0, \delta > 0$.}
	\textbf{Pre-compute} $b_{j} = \sum_{i} d_{i}h_{ij}$, $\forall j$. \\
	\For {$n = 0, 1, 2, ...$}{
		$q_{i}^{(n)} = I_{0,i} \exp(-\sum_{j}h_{ij}x^{(n)}_{j})$, $\forall i$. \\
		${b}^{(n)}_{j} = \sum_{i} q_{i}^{(n)} h_{ij}$, $\forall j$. \\
		$x_{j}^{(n+1)} = \Big[x^{(n)}_{j} - \frac{1}{L_{\phi}}  \Big( (b_{j} - {b}^{(n)}_{j}) + \lambda \frac{\partial \beta(\boldsymbol{x}^{(n)})}{\partial x^{(n)}_{j}} \Big) \Big]_{+}    $, $\forall j$.
	}
\end{algorithm}

The acceleration methods have been investigated using two real-data scans acquired on a SureScan\textsuperscript{TM} \textit{x}1000 Explosive Detection System. For all the methods presented, the regularization parameters $\lambda = 15000$ and $\delta = 0.001$ were used. All reconstructed image volumes consist of $240$ rows, $400$ columns, and a varying number of slices along the z-axis. In order to compare the rates of convergence, we look at the normalized function errors $(\Phi(\boldsymbol{x}^{(n)}) - \Phi(\boldsymbol{x}^{*}))/ \|\Phi(\boldsymbol{x}^{*}) \|$, where the objective function value at the optimum was computationally found by running the algorithms for many more iterations. We examine the performance of our algorithms with a relatively low number of subsets, $8$, and a large number of subsets, $64$.

We compare the proposed algorithms SA-JS and OSA-JS with
\begin{itemize}
	\item the full variant of Jensen surrogates, Full-JS,
	\item the momentum and varying step-size variant, Fast-JS\cite{degirmenci2016thesis},
	\item the ordered subsets with Jensen Surrogates method, OS-JS,
	\item the Stochastic Average Gradient method,\cite{roux2012stochastic} SA-GD, and
	\item the ordered subsets with gradient descent method, OS-GD.
\end{itemize}

\begin{figure}
	\captionsetup{justification=centering}
	\centering
	\begin{subfigure}{.3\textwidth}
		\includegraphics[width=1.0\linewidth]{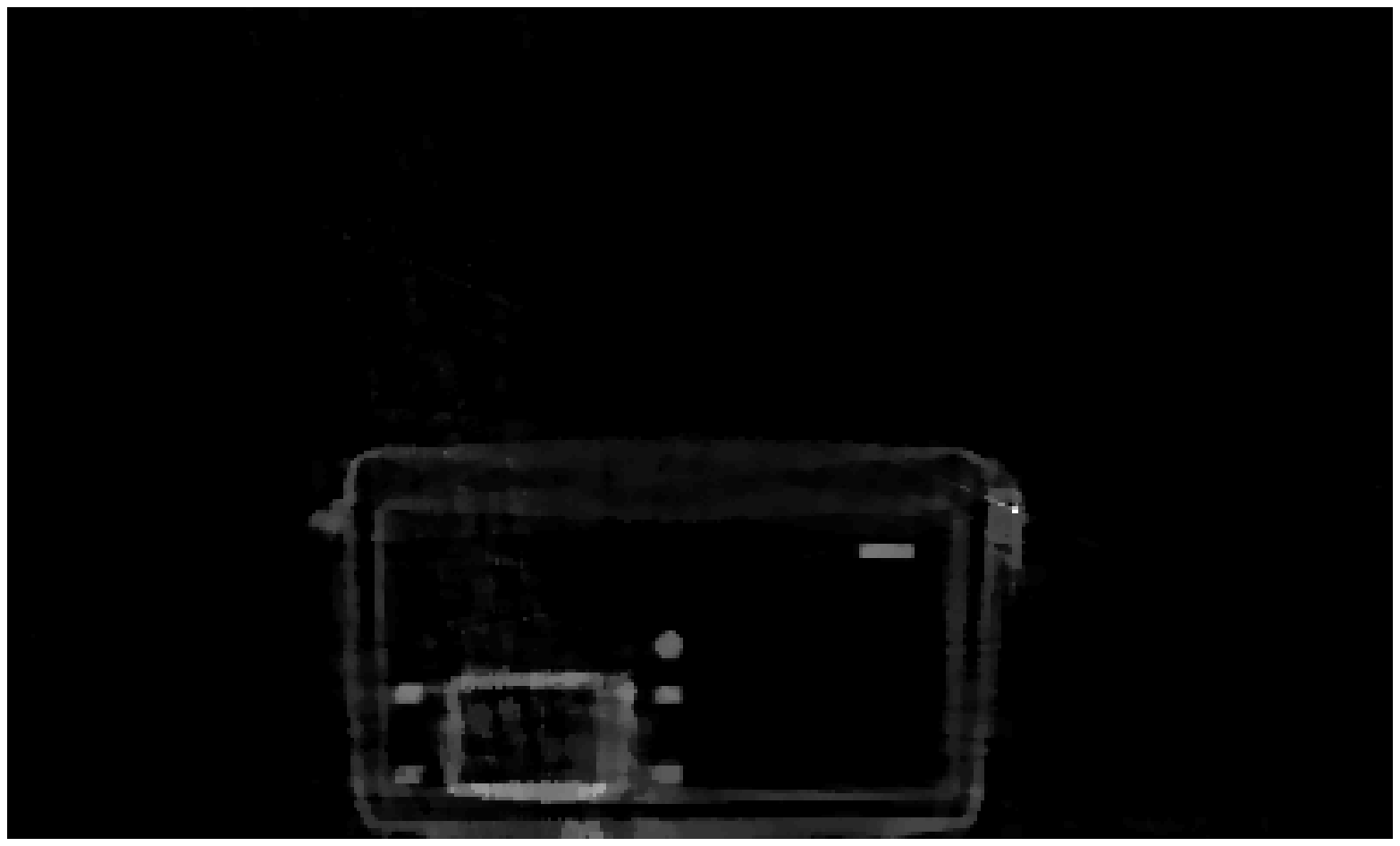}
	\end{subfigure}
	\begin{subfigure}{.3\textwidth}
		\includegraphics[width=1.0\linewidth]{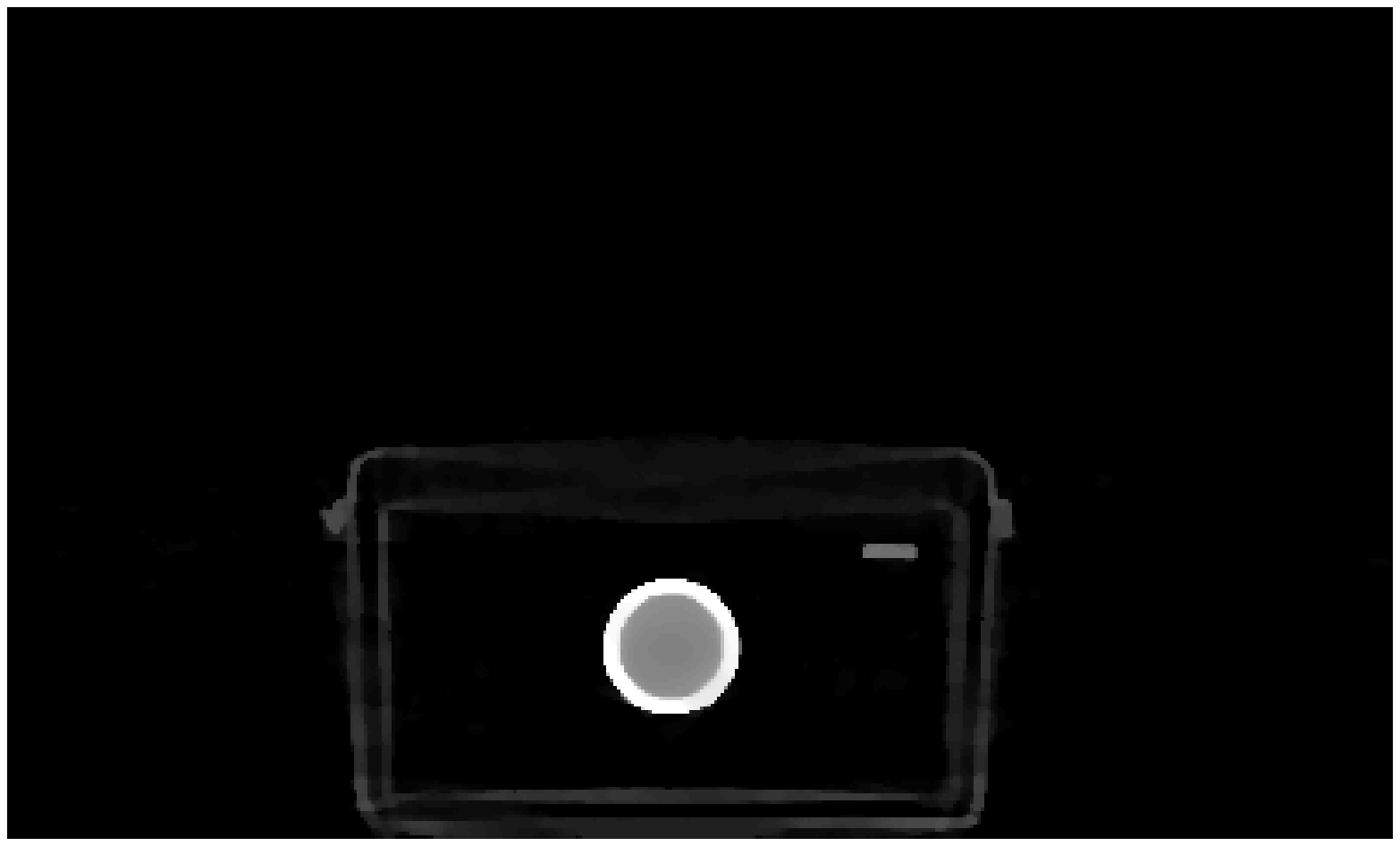}
	\end{subfigure}
	\begin{subfigure}{.3\textwidth}
		\includegraphics[width=1.0\linewidth]{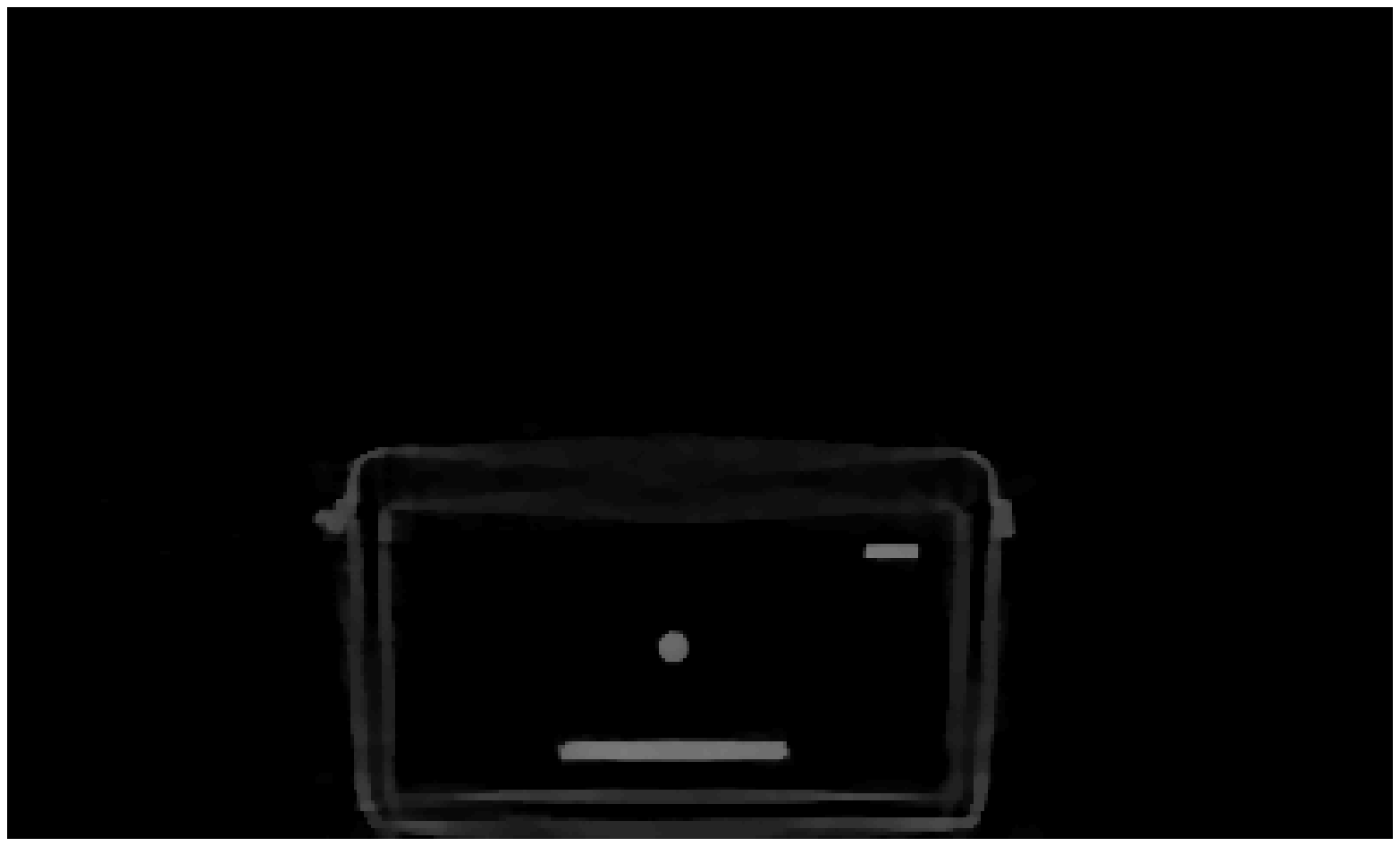}
	\end{subfigure}
	\caption{Different axial views for the reconstructed volume for Bag \#1.}
	\label{fig1}
\end{figure}

\begin{figure}
	\captionsetup{justification=centering}
	\centering
	\begin{subfigure}{.3\textwidth}
		\includegraphics[width=1.0\linewidth]{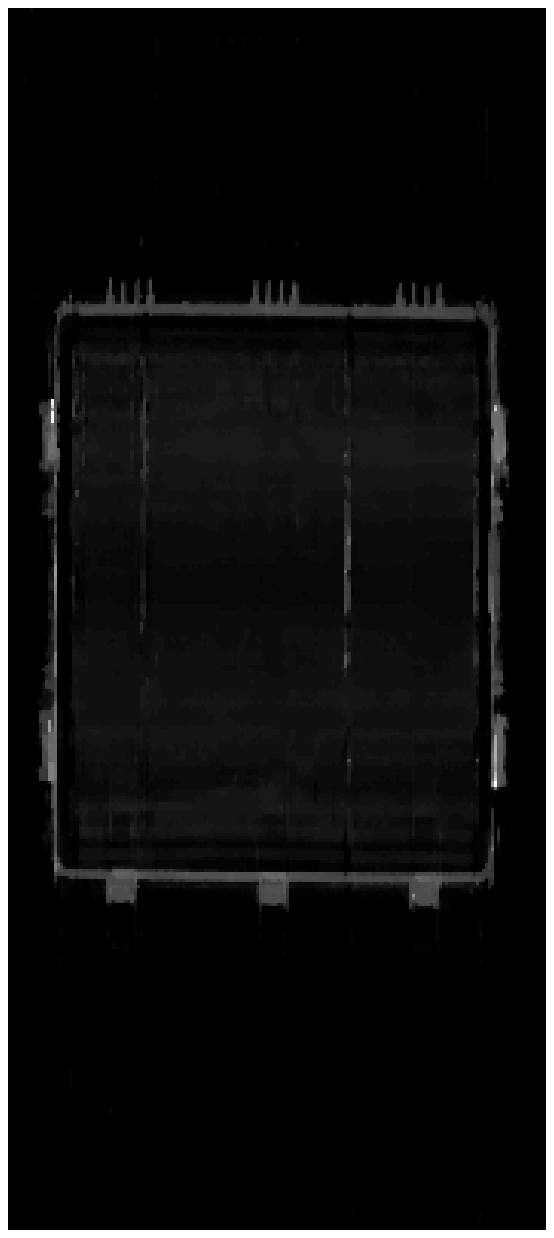}
	\end{subfigure}
	\begin{subfigure}{.3\textwidth}
		\includegraphics[width=1.0\linewidth]{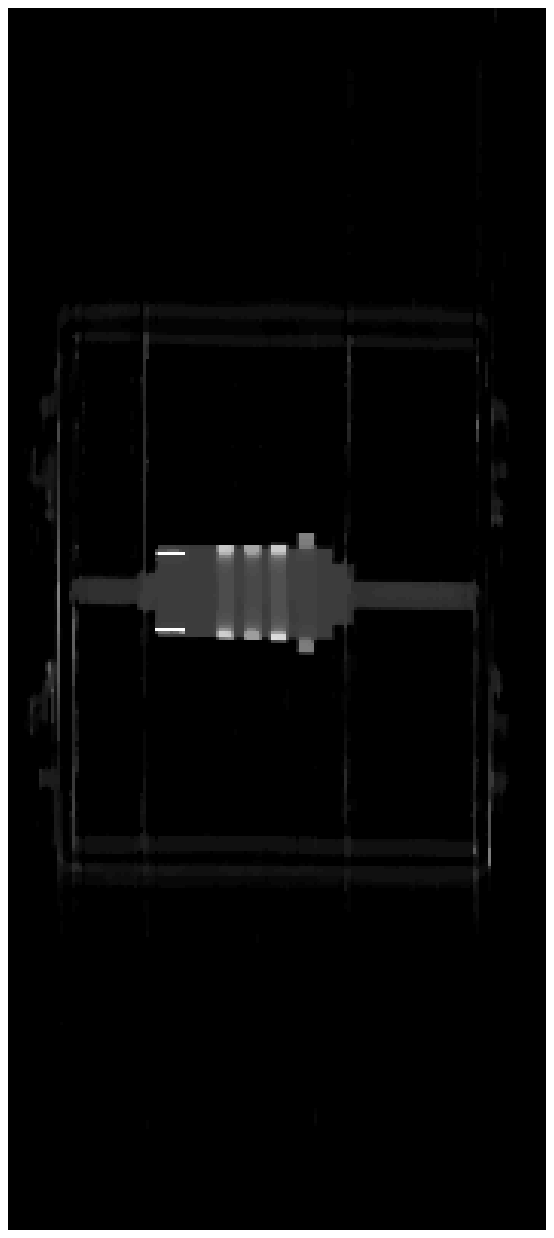}
	\end{subfigure}
	\begin{subfigure}{.3\textwidth}
		\includegraphics[width=1.0\linewidth]{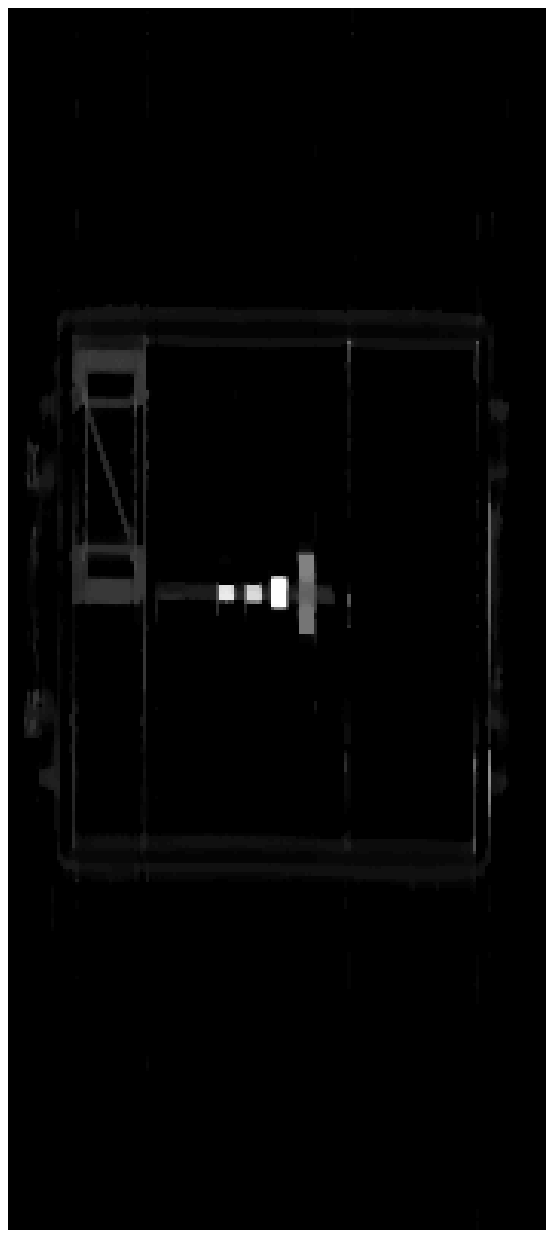}
	\end{subfigure}
	\caption{Different sagittal views for the reconstructed volume for Bag \#1.}
	\label{fig2}
\end{figure}

\begin{figure}
	\captionsetup{justification=centering}
	\centering
	\begin{subfigure}{.3\textwidth}
		\includegraphics[width=1.0\linewidth]{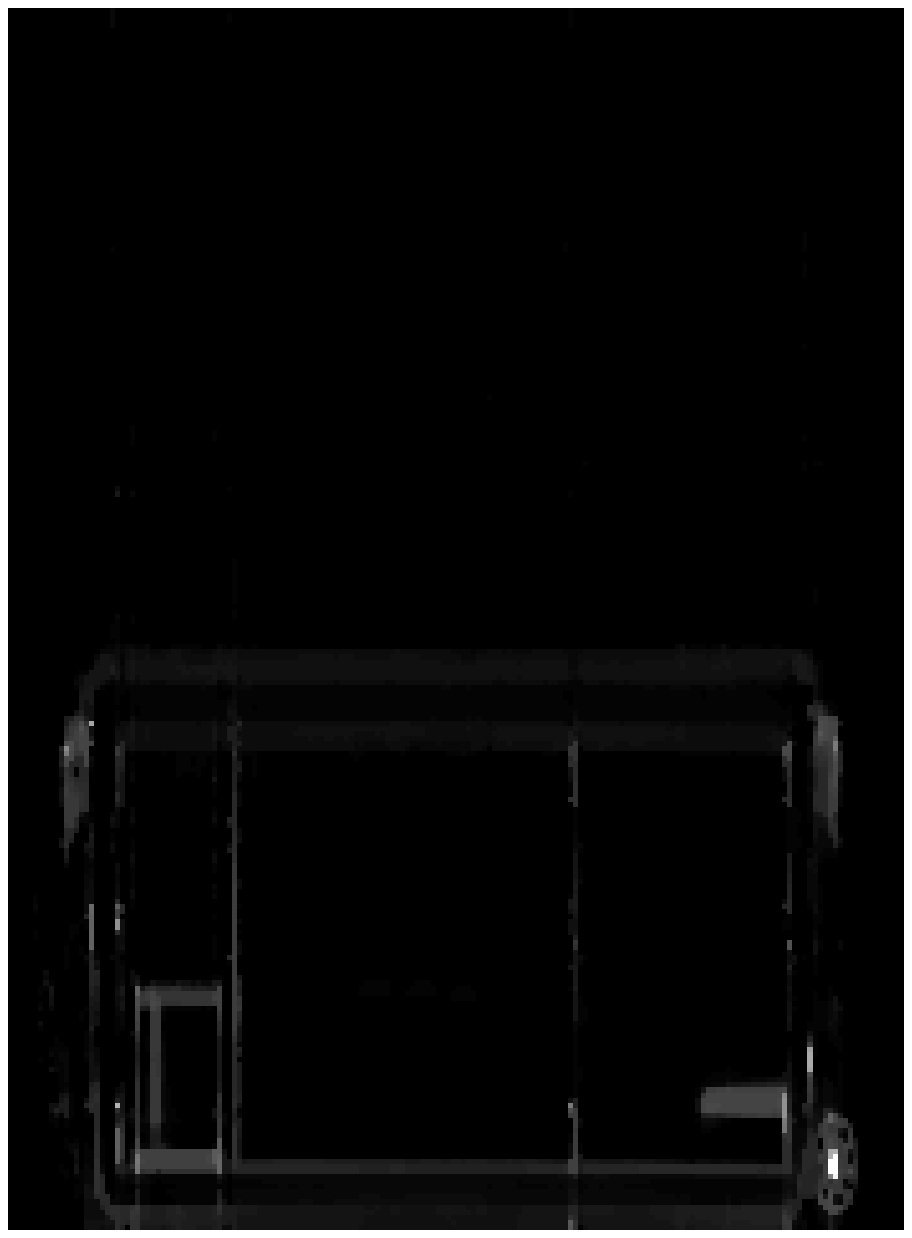}
	\end{subfigure}
	\begin{subfigure}{.3\textwidth}
		\includegraphics[width=1.0\linewidth]{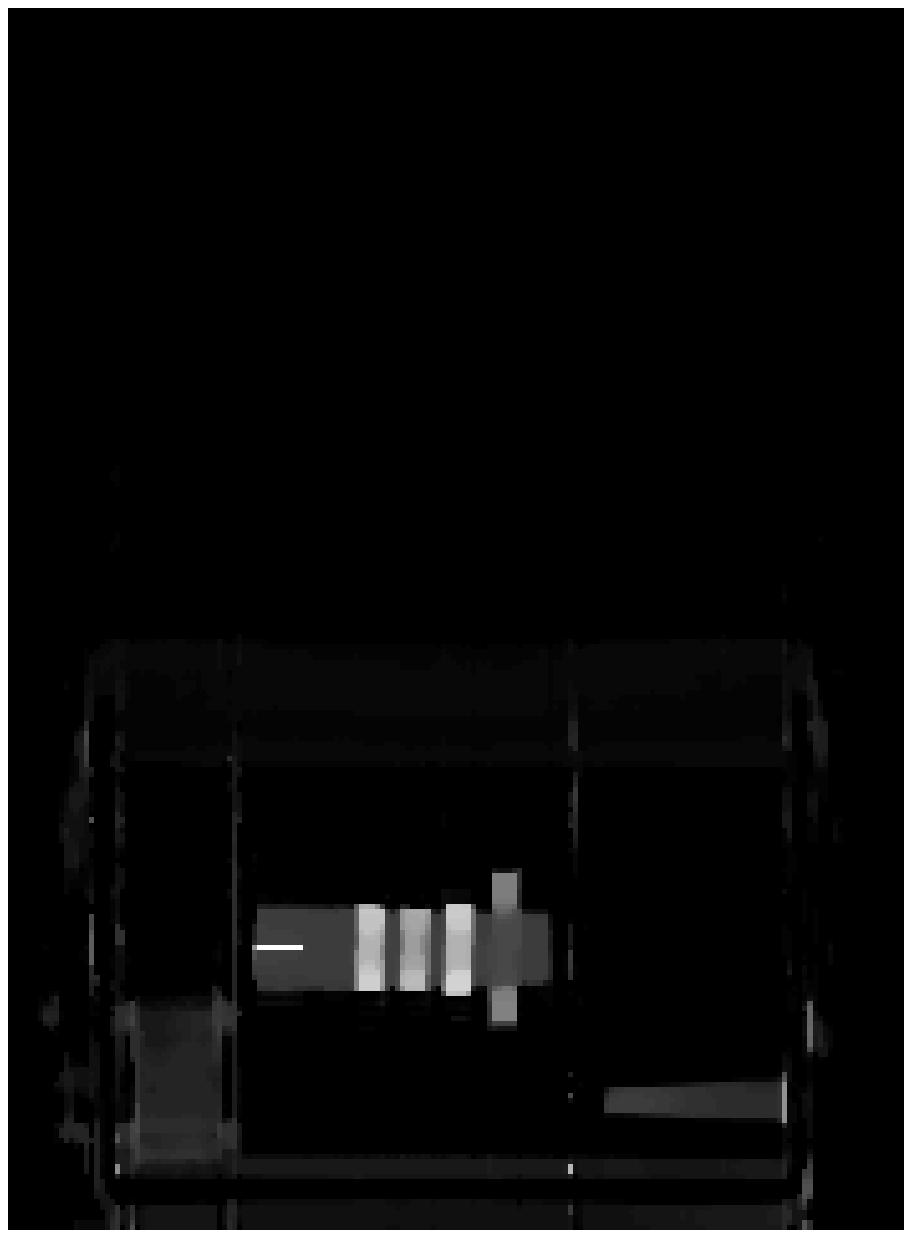}
	\end{subfigure}
	\begin{subfigure}{.3\textwidth}
		\includegraphics[width=1.0\linewidth]{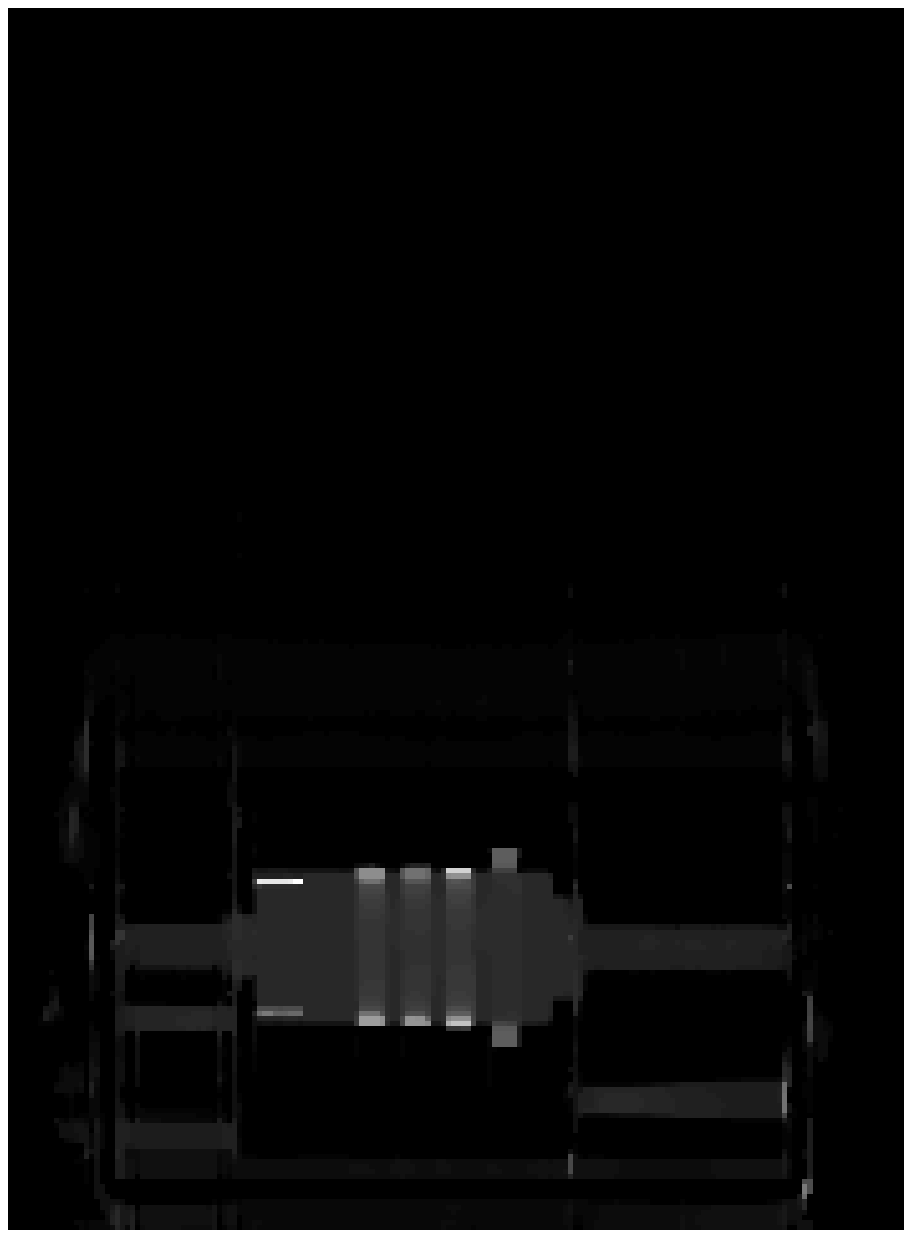}
	\end{subfigure}
	\caption{Different coronal views for the reconstructed volume for Bag \#1.}
	\label{fig3}
\end{figure}

\begin{figure} 
	\captionsetup{justification=centering}
	\centering
	\begin{subfigure}{.45\textwidth}
		\includegraphics[width=1.0\linewidth]{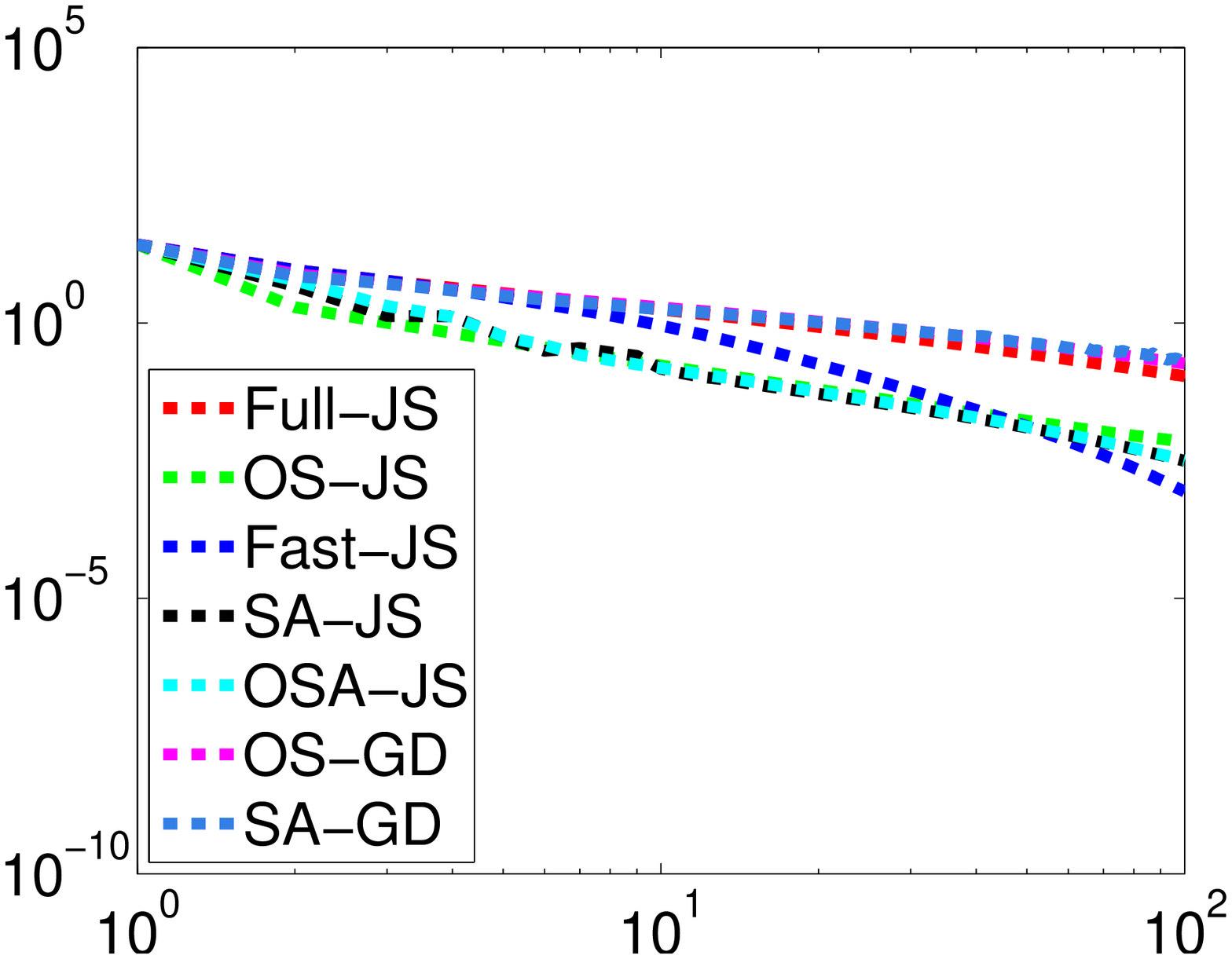}
		\caption{8 Subsets}
	\end{subfigure}
	\begin{subfigure}{.45\textwidth}
		\includegraphics[width=1.0\linewidth]{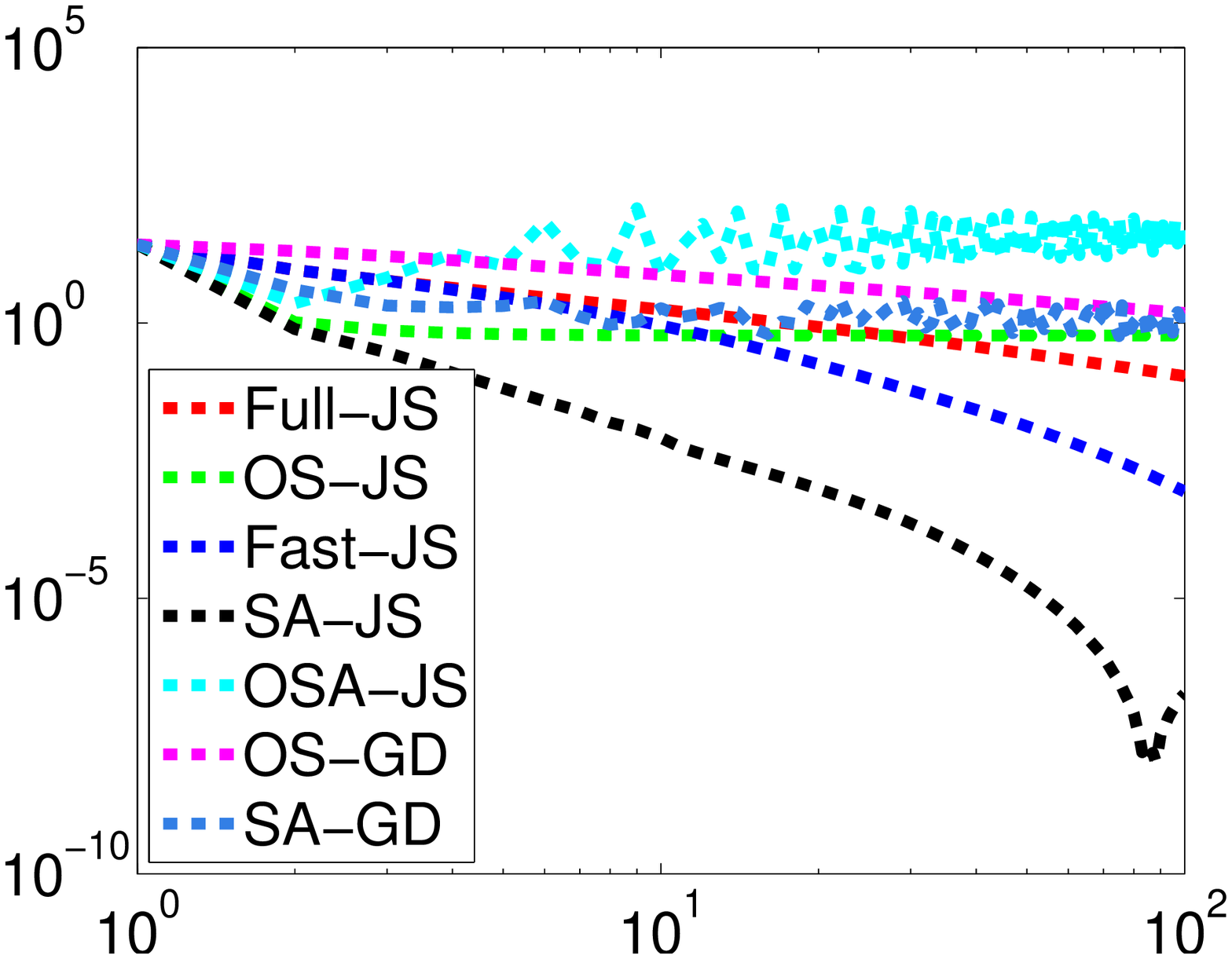}
		\caption{64 Subsets}
	\end{subfigure}
	\caption{Normalized function errors vs. number of effective data passes for Bag \#1. (a) shows the results for $8$ subsets, while (b) shows the results for $64$ subsets.}
	\label{fig4}
\end{figure}

\begin{figure} 
	\captionsetup{justification=centering}
	\centering
	\begin{subfigure}{.45\textwidth}
		\includegraphics[width=1.0\linewidth]{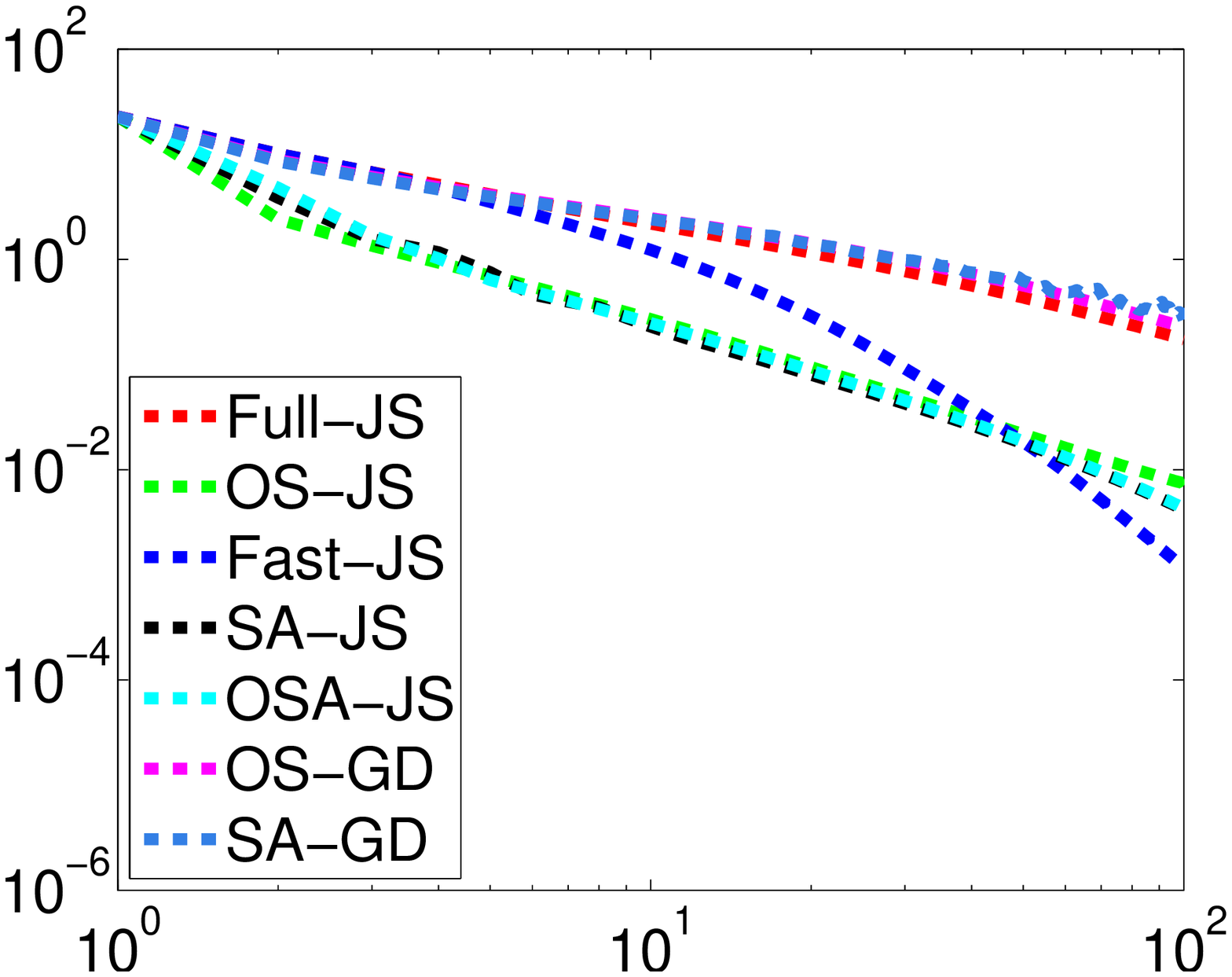}
		\caption{8 Subsets}
	\end{subfigure}
	\begin{subfigure}{.45\textwidth}
		\includegraphics[width=1.0\linewidth]{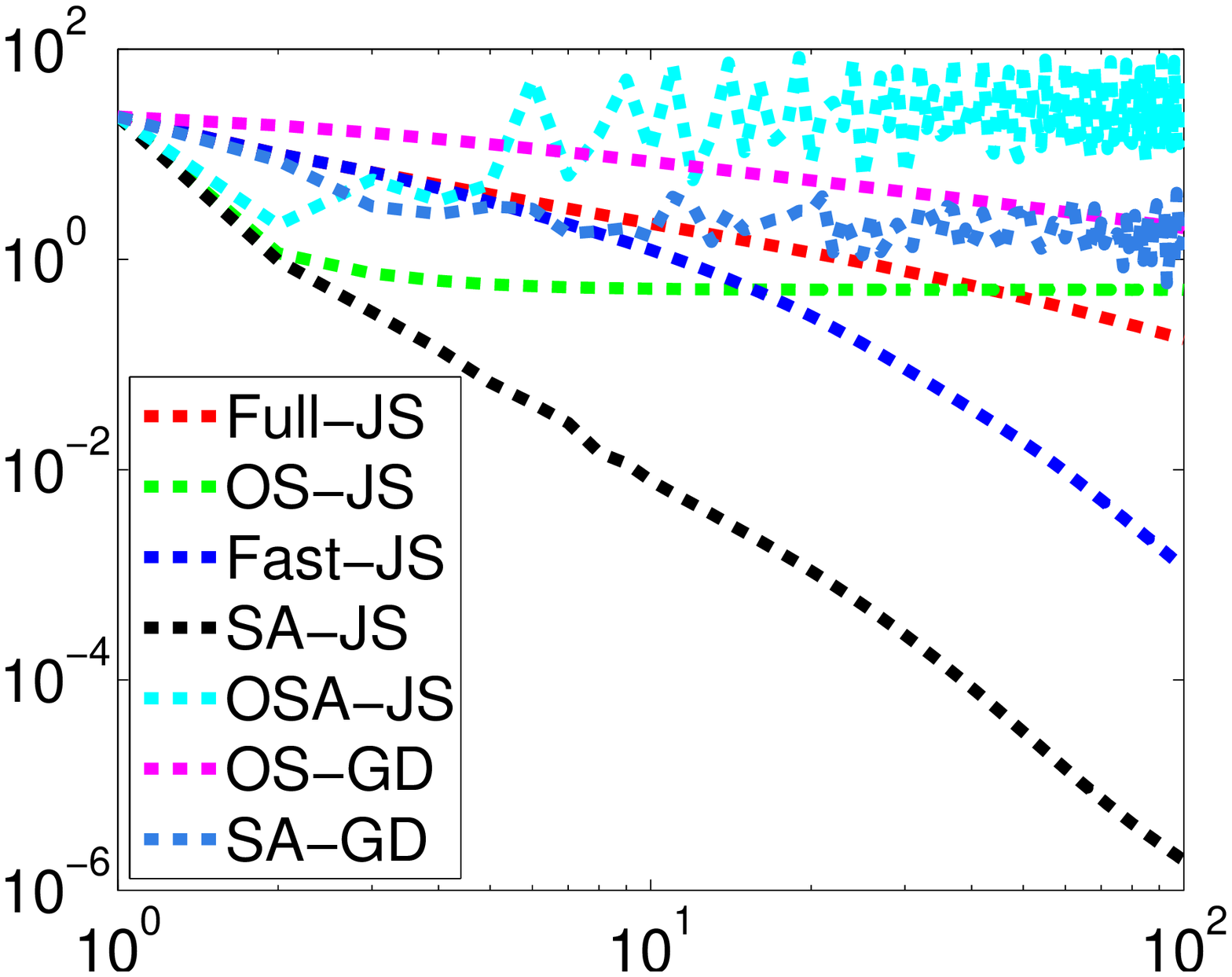}
		\caption{64 Subsets}
	\end{subfigure}
	\caption{Normalized function errors vs. number of effective data passes for Bag \#2. (a) shows the results for $8$ subsets, while (b) shows the results for $64$ subsets.}
	\label{fig8}
\end{figure}

Figures~\ref{fig1}-\ref{fig3} show three slices each of axial, sagittal, and coronal views of the reconstructed volume for Bag \#1. Figure~\ref{fig4} presents the normalized function errors for different algorithms vs. the number of effective data passes for Bag \#1 for two subset settings while Figure~\ref{fig8} presents the errors for Bag \#2. As seen in the figures, for $8$ subsets, the Fast-JS algorithm slightly outperforms the SA-JS and OSA-JS algorithms proposed. For $64$ subsets, we see that SA-JS algorithm we proposed outperforms the other competing algorithms.

It is also important to note that the proposed Jensen surrogate algorithms outperform their gradient descent equivalents for all cases. Another observation is that the cyclic choice we proposed for the average scheme we proposed works competitively for the $8$ subsets case, but performs worse than its stochastic counterpart for $64$ subsets.

\section{Conclusion}

In this paper, we investigated incremental methods and proposed one stochastic and one deterministic algorithm that use Jensen surrogates and compared these methods with other competing algorithms using data collected from two bags by a baggage scanner. The Stochastic Average variant (SA-JS) we proposed performs competitively for the small number of subsets and outperforms the other competing algorithms for the large number of subsets in terms of number of effective data passes. This is an indicator that it should perform better as the number of subsets increases. However, implementation of the larger number of subsets requires additional storage and more image updates per effective data pass, which both add computational costs. Thus, we plan to conduct a thorough time analysis of the proposed algorithms as well as to explore the ``hybrid" variants.

\bibliography{dissertation_bibliography_111115}   
\bibliographystyle{spiebib}   
\end{document}